\documentclass[12pt]{article}
\usepackage{amssymb,latexsym,amsmath}
\def\ZZ{{\mathbb Z}}
\def\QQ{{\mathbb Q}}
\def\RR{{\mathbb R}}
\def\CC{{\mathbb C}}

\newtheorem{formula}{}[section]
\newtheorem{proposition}[formula]{Proposition}
\newtheorem{definition}[formula]{\indent Definition}
\newtheorem{corollary}[formula]{\indent Corollary}
\newtheorem{remark}[formula]{\indent Remark}
\newtheorem{lemma}[formula]{\indent Lemma}
\newtheorem{theorem}[formula]{\indent Theorem}

\newtheorem{example}[formula]{Example}

\def\thrm{\begin{theorem}}
\def\thrml#1{\begin{theorem}\label{#1}}
\def\ethrm{\end{theorem}}
\def\rmrk{\begin{remark}}
\def\rmrkl#1{\begin{remark}\label{#1}}
\def\ermrk{\end{remark}}
\def\dfntn{\begin{definition}}
\def\dfntnl#1{\begin{definition}\label{#1}}
\def\edfntn{\end{definition}}
\def\nmrt{\begin{enumerate}}
\def\enmrt{\end{enumerate}}

\def\qtn{\begin{equation}}
\def\qtnl#1{\begin{equation}\label{#1}}
\def\eqtn{\end{equation}}
\def\lmm{\begin{lemma}}
\def\lmml#1{\begin{lemma}\label{#1}}
\def\elmm{\end{lemma}}
\def\crllr{\begin{corollary}}
\def\crllrl#1{\begin{corollary}\label{#1}}
\def\ecrllr{\end{corollary}}

\begin{document}
\title{}
\date{}
\maketitle
%\nopagenumbers
%\begin{titlepage}
%\title
\vspace{-0,1cm} \centerline{\bf Min-plus Bruhat decomposition}
%\centerline{\bf AND RELATIVE KOLMOGOROV COMPLEXITY}
\vspace{7mm}
\author{
\centerline{Dima Grigoriev}
%\\[-1pt]
\vspace{3mm}
%\centerline{$^1$ St.Petersbourg University,
% Universitetskaya nab., 7/9,}
%\centerline{ St.Petersbourg,
% 199164, RUSSIA}
%\vspace{3mm}
%Dima Grigoriev \\[-1pt]
\centerline{CNRS, Math\'ematique, Universit\'e de Lille, Villeneuve
d'Ascq, 59655, France} \vspace{1mm} \centerline{e-mail:\
dmitry.grigoryev@univ-lille.fr } \vspace{1mm}
\centerline{URL:\ http://en.wikipedia.org/wiki/Dima\_Grigoriev} }
%\date{}
%\maketitle

\begin{abstract}
We study min-plus Bruhat cells $B_{trop}^{(n)} \odot p\odot B_{trop}^{(n)}$, where $B_{trop}^{(n)}$ denotes the set of all upper-triangular matrices over $\RR \cup \{\infty\}$, and $p$ is a permutation. It is proved that $\dim (B_{trop}^{(n)} \odot p \odot B_{trop}^{(n)})= \frac{n(n+1)}{2} + l(p)$, where $l(p)$ is the length of $p$ in Coxeter generators. Also we provide an explicit description of the cell $B_{trop}^{(n)}(\RR) \odot p_0\odot B_{trop}^{(n)}(\RR)$ (up to a subset of dimension less than $n^2$) for the highest in Bruhat order permutation $p_0$. Finally, it is shown that for a tropically non-singular matrix $A$ for at most one permutation $p_A$ it can hold $A\in B_{trop}^{(n)}\odot p_A\odot B_{trop}^{(n)}$. 
\end{abstract}

{\bf keywords}: min-plus Bruhat decomposition, dimension of min-plus Bruhat cell


\section{Introduction}
We consider min-plus algebra over $\RR_{\infty}:= \RR \cup \{\infty\}$ with operations $\oplus:= \min,\ \odot:= +$ (one can find basic concepts of min-plus algebra in \cite{MS}, \cite{J}). Our goal is to study a min-plus analogue of Bruhat decomposition. Recall that the latter in classical algebra states that every  $n\times n$ matrix $M$ belongs to Bruhat cell $M\in B^{(n)}p B^{(n)}$ for a  permutation $p\in S_n$, where $B^{(n)}$ denotes the set of  upper-triangular matrices. When $M$ is a non-singular matrix, the permutation is unique. We identify $p$ with the matrix $p$ whose $(i, p(i))$-entries equal $0$ for $1\le i\le n$, and all other entries equal $\infty$. Chevalley's theorem \cite{St}, \cite{F} states that the closure
$$\overline{B^{(n)}pB^{(n)}}=\bigcup_{p_1\preceq p} B^{(n)}p_1 B^{(n)},$$
\noindent where  $\prec$ denotes Bruhat order on $S_n$. Originally, Chevalley's theorem provided this equality after restricting to non-singular matrices, the extension to arbitrary matrices one can find in \cite{G81}, \cite{G82}.

Denote by $B_{trop}^{(n)}$ the set of all $n\times n$ upper-triangular matrices over $\RR_{\infty}$. Observe that a min-plus Bruhat cell $B_{trop}^{(n)}\odot p\odot B_{trop}^{(n)}$ is a min-plus (or tropical) prevariety \cite{MS}, \cite{J}.
If a min-plus matrix $A=b\odot p\odot c$ for some $p\in S_n$ and $b,c\in B_{trop}^{(n)}$, we say that $A$ admits a min-plus Bruhat decomposition. Note that not every min-plus matrix admits a min-plus Bruhat decomposition. For example, one can easily see that for $n=2$ the min-plus Bruhat cell $B_{trop}^{(2)}\odot (12) \odot B_{trop}^{(2)}$ consists exactly of matrices $(a_{i,j})_{1\le i,j\le 2}$ such that $a_{1,1}+a_{2,2}\ge a_{1,2}+a_{2,1}$.

For a permutation $p\in S_n$ by $l(p)$ denote its length in Coxeter generators $(i\ i+1),\ 1\le i<n$ \cite{St} (which equals the number of inversions in $p$). In section~\ref{one} we prove the following

\begin{theorem}\label{dimension}
$$\dim(B_{trop}^{(n)}\odot p\odot B_{trop}^{(n)})=\frac{n(n+1)}{2} +l(p).$$
\end{theorem}

In section~\ref{two} we study the highest min-plus Bruhat cell $B_{trop}^{(n)}\odot p_0\odot B_{trop}^{(n)}$ where $p_0=(1\ n) (2\ n-1)\cdots$ is the (unique) permutation with the largest possible length $l(p_0)=\frac{n(n-1)}{2}$, being the highest in Bruhat order: $p\prec p_0$ for every permutation $S_n\ni p\neq p_0$. Theorem~\ref{dimension} implies that $\dim(B_{trop}^{(n)}  \odot p_0 \odot B_{trop}^{(n)})=n^2$.
Let $A=(a_{i,j})_{1\le i,j\le n},\ a_{i,j}\in \RR$ be a finite matrix. Denote 
$$\Delta_{q,s}^{k,m}:= a_{q,k}+a_{s,m}-a_{q,m}-a_{s,k},\ 1\le q<s\le n,\ 1\le k<m\le n.$$
\noindent Note that $\Delta_{q,s}^{k,m}=0$ means that the corresponding $2\times 2$ tropical submatrix of $A$ is tropically singular (see e.g. \cite{AGG}). The following theorem provides an explicit description of the cell $B_{trop}^{(n)}  \odot p_0 \odot B_{trop}^{(n)} \cap \RR^{n^2}$ up to a subset of dimension less than $n^2$. Note that not finite matrices constitute a tropical prevariety of dimension less than $n^2$ as well.

\begin{theorem}\label{highest}
i) For an $n\times n$ finite matrix $A\in \RR^{n^2}$ if it holds
%it holds $A\in B_{trop}^{(n)}  \odot p_0 \odot B_{trop}^{(n)}$ iff 
\begin{equation}\label{2}
\Delta_{q,s}^{k,m}\ge 0,\ 1\le q<s\le n,\ 1\le k<m\le n,\ k+s=n+1.
\end{equation}
then $A$ admits the
 following canonical decomposition:
\begin{equation}\label{1}
A=(b_{i,j})_{1\le i,j\le n} \odot p_0\odot (c_{i,j})_{1\le i,j\le n}\in B_{trop}^{(n)}\odot p_0\odot B_{trop}^{(n)},
\end{equation} 
where $b_{i,j}:=a_{i,n+1-j}-a_{j,n+1-j},\ c_{i,j}:=a_{n+1-i,j},\ 1\le i\le j\le n$ and $b_{i,j}=c_{i,j}=\infty,\ 1\le j<i\le n$.

ii) Conversely, if $A\in B_{trop}^{(n)}\odot p_0\odot B_{trop}^{(n)}$ then either
 all the inequalities (\ref{2}) are strict:
\begin{equation}\label{6}
\Delta_{q,s}^{k,m}> 0,\ 1\le q<s\le n,\ 1\le k<m\le n,\ k+s=n+1
\end{equation}
or $\Delta_{q,s}^{k,m}=0$ for some $1\le q<s\le n,\ 1\le k<m\le n,\ k+s=n+1$.

iii) If (\ref{6}) holds then
the canonical decomposition (\ref{1}) is unique. Namely, if
$$A=(b_{i,j}')_{1\le i,j\le n}\odot p_0 \odot (c_{i,j}')_{1\le i,j\le n},\ b_{i,i}'=0,\ b_{i,j}'=c_{i,j}'=\infty,\ 1\le j<i\le n$$
\noindent then $b_{i,j}'=b_{i,j},\ c_{i,j}'=c_{i,j},\ 1\le i,j\le n$ where $b_{i,j},\ c_{i,j}$ are from (\ref{1}).
\end{theorem} 

In section~\ref{three} for a tropically non-singular (see e.g. \cite{AGG}) matrix $A$ we show the uniqueness of a permutation $p_A$ (provided that it does exist) such that $A\in B_{trop}^{(n)}\odot p_A\odot  B_{trop}^{(n)}$. This is a min-plus analogue of the assertion that $B_1^{(n)}p_1 B_1^{(n)} \cap B_1^{(n)}p_2 B_1^{(n)} =\emptyset,\ p_1\neq p_2\in S_n$.
 
We mention that in \cite{N} a min-plus version of Bruhat decomposition was considered for matrices over the semi-ring $\RR_{\infty} \cup G$ where $G$ is a set of "ghosts" isomorphic to $\RR$ and "merged" with $\RR$, so that between any pair of "adjacent" real numbers a "ghost" is inserted.

 \section{Dimension of a min-plus Bruhat cell}\label{one}

{\bf Proof of Theorem~\ref{dimension}}. First we prove the upper bound on the dimension
$$\dim (B_{trop}^{(n)} \odot p \odot B_{trop}^{(n)}) \le \frac{n(n+1)}{2}+l(p).$$
\noindent Now we denote by $B^{(n)}$ the set of all upper-triangular $n\times n$ matrices over the field $\CC ((t^{1/\infty}))$ of Newton-Puiseux series of the form $\sum_{i_0\le i\in \ZZ} \alpha_i t^{i/q},\ \alpha_i \in \CC,\ q\in \ZZ_{\ge 1}$.

We claim that 
$$B_{trop}^{(n)} \odot p \odot B_{trop}^{(n)} \subset Trop(B^{(n)}pB^{(n)}),$$
\noindent where $Trop$ denotes the tropicalization \cite{MS}, \cite{J}. Indeed, for any matrix 
$$(b_{i,j})_{1\le i,j\le n} \odot p\odot (c_{i,j})_{1\le i,j\le n} \in B_{trop}^{(n)} \odot p \odot B_{trop}^{(n)}\cap (\QQ_{\infty})^{n^2}$$
\noindent with $b_{i,j},\  c_{i,j}\in \QQ_{\infty}$ consider (classical) matrices  $(t^{b_{i,j}})_{1\le i,j\le n},\ (t^{c_{i,j}})_{1\le i,j\le n}$.
%, where $\beta_{i,j},\  \gamma_{i,j}\in \CC$ are $\QQ$-linearly independent. 
As usually, we agree that $t^{\infty}:=0$. Observe that if for a matrix $A$ over the semi-ring $\QQ_{\infty}$ it holds
\begin{equation}\label{5}
A=(b_{i,j})_{1\le i,j\le n} \odot p\odot (c_{i,j})_{1\le i,j\le n}
\end{equation}
for some $b_{i,j},\ c_{i,j}\in \RR_{\infty}$ then there exists a solution of (\ref{5}) for suitable $b_{i,j},\ c_{i,j}\in \QQ_{\infty}$ since (\ref{5}) is a system of (classical) linear inequalities. Thus, we obtain the equality
$$Trop((t^{b_{i,j}})_{1\le i,j\le n}\cdot p\cdot (t^{c_{i,j}})_{1\le i,j\le n})=(b_{i,j})_{1\le i,j\le n} \odot p\odot (c_{i,j})_{1\le i,j\le n},$$
\noindent which justifies the claim.

On the other hand, it holds
$$\dim(Trop(B^{(n)}pB^{(n)}))= \dim(B^{(n)}pB^{(n)})=\frac{n(n+1)}{2}+l(p) .$$
\noindent For the first equality see \cite{MS}, \cite{J}. For the second equality see \cite{St}, \cite{F} (one can deduce it by induction on Bruhat order from Chevalley's theorem, see also \cite{G81}, \cite{G82}). Together with the claim above this implies the upper bound on the dimension. \vspace{1mm} 

Now we proceed to the lower bound. We consider the equality
$$(a_{i,j})_{1\le i,j\le n}=(b_{i,j})_{1\le i,j\le n}\odot p\odot (c_{i,j})_{1\le i,j\le n}\in B_{trop}^{(n)} \odot p \odot B_{trop}^{(n)}.$$
\noindent It holds 
$$a_{i,j}=\min_{i\le k\le n,\ 1\le p(k)\le j} \{b_{i,k}+c_{p(k),j}\},\ 1\le i,j\le n.$$
We impose equalities $b_{i,i}=0,\ b_{i,k}=\infty,\ 1\le i<k\le n,\ p(i)<p(k)$. 

Let us construct a min-plus prevariety of entries $a_{i,j}$ by recursion on $n-i$ (simultaneously we construct matrices $(b_{i,j})_{1\le i,j\le n},\ (c_{i,j})_{1\le i,j\le n}$). For the base of recursion $i=n$ we have
$$a_{n,j}=c_{p(n),j},\ p(n)\le j\le n,\ a_{n,q}=\infty,\ 1\le q<p(n),$$
\noindent which determines $c_{p(n),j},\ p(n)\le j\le n$ for arbitrary $a_{n,j},\ p(n)\le j\le n$.

For the recursive hypothesis assume that for some $1\le i<n$ for all $n\ge m>i$ entries
$a_{m,r},\ b_{m,q},\ c_{p(m),j},\ 1\le r,q,j\le n$ are already determined, and $b_{m,q}<\infty,\ m\le q,\ p(m)\ge p(q),\ c_{p(m),j}<\infty,\ p(m)\le j$. Let $p(m_1)<\cdots < p(m_s)<p(i)$ for $m_1,\dots,m_s>i$, and $p(m)>p(i)$ for any $m>i,\ m\notin \{m_1,\dots,m_s\}$. For $q<p(m_1)$ it holds $a_{i,q}=\infty$. For $p(m_1)\le q<p(m_2)$ it holds $a_{i,q}=b_{i,m_1}+c_{p(m_1),q}$. Therefore, for an arbitrary value of $a_{i,p(m_1)}<\infty$ we determine $b_{i,m_1}:=a_{i,p(m_1)}-c_{p(m_1),p(m_1)}$. This determines also $a_{i,q},\ p(m_1)<q<p(m_2)$. After that for $p(m_2)\le q<p(m_3)$ it holds $a_{i,q}=\min\{b_{i,m_1}+c_{p(m_1),q},\ b_{i,m_2}+c_{p(m_2),q}\}$. For an arbitrary value of $a_{i,p(m_2)}<b_{i,m_1}+c_{p(m_1),p(m_2)}$ we have $a_{i,p(m_2)}=b_{i,m_2}+c_{p(m_2),p(m_2)}$, which determines $b_{i,m_2}:=a_{i,p(m_2)}-c_{p(m_2),p(m_2)}$. This also determines $a_{i,q},\ p(m_2)<q<p(m_3)$. Continuing in this way, we take arbitrary values of $a_{i,p(m_1)}, a_{i,p(m_2)},\dots, a_{i, p(m_s)}$ (in suitable semi-infinite intervals), which determine $b_{i,m_1},\dots,b_{i,m_s}$, respectively, and $a_{i,q},\ 1\le q<p(i)$.

Then for every $p(i)\le q\le n$ it holds
$$a_{i,q}=\min\{b_{i,m_1}+c_{p(m_1),q},\dots, b_{i,m_s}+c_{p(m_s),q},\ c_{p(i),q}\}.$$
\noindent For an arbitrary value of $a_{i,q}< \min\{b_{i,m_1}+c_{p(m_1),q},\dots, b_{i,m_s}+c_{p(m_s),q}\}$ it holds $a_{i,q}=c_{p(i),q}$, which determines $c_{p(i),q}$. This completes the recursive step.

Summarizing the recursive construction, the values of $a_{i,q},\ 1\le i\le n,\ p(i)\le q\le n$ were taken arbitrarily (each of them from a semi-infinite interval), in addition, the values of $a_{i,p(m)},\ m>i,\ p(m)<p(i)$ also were taken arbitrarily (each of them from a semi-infinite interval). This provides the lower bound
$$\dim(B_{trop}^{(n)}\odot p\odot  B_{trop}^{(n)})\ge \frac{n(n+1)}{2} + l(p).$$ $\Box$

\section{The highest min-plus Bruhat cell}\label{two}

{\bf Proof of Theorem~\ref{highest}}.  i) First assume (\ref{2}). We prove (\ref{1}) by induction on $n-i$. For the base of induction $i=n$ we take $c_{1,j}:=a_{n,j},\ 1\le j\le n$.

For the inductive step assume that (\ref{1}) and the canonical decomposition are already established for all rows $r>i$. It holds $a_{i,1}=(a_{i,1}-a_{n,1})+a_{n,1}=b_{i,n}+c_{1,1}$. After that we have 
$$a_{i,2}=\min\{(a_{i,1}-a_{n,1})+a_{n,2},\ (a_{i,2}-a_{n-1,2})+a_{n-1,2}\}=\min\{b_{i,n}+c_{1,2},\ b_{i,n-1}+c_{2,2}\},$$
\noindent which is true since $\Delta_{i,n}^{1,2}\ge 0$.  Similarly, for $j\le n-i+1$ we get
$$\min\{b_{i,n}+c_{1,j},\dots, b_{i,n-j+2}+c_{j-1,j},\ b_{i,n-j+1}+c_{j,j}\}=$$
\begin{equation}\label{3}
\min\{a_{i,1}-a_{n,1}+a_{n,j},\dots, a_{i,j-1}-a_{n-j+2,j-1}+a_{n-j+2,j},\ a_{i,j}-a_{n-j+1,j}+a_{n-j+1,j}\}=a_{i,j},
\end{equation}
which is valid because $\Delta_{i,n}^{1,j},\ \Delta_{i,n-1}^{2,j},\dots, \Delta_{i,n-j+2}^{j-1,j}\ge 0$.

For $n-i+2\le j\le n$ we have 
$$\min\{b_{i,n}+c_{1,j},\dots, b_{i,i+1}+c_{n-i,j},\ b_{i,i}+c_{n-i+1,j}\}=$$
\begin{equation}\label{4}
\min\{a_{i,1}-a_{n,1}+a_{n,j},\dots, a_{i,n-i}-a_{i+1,n-i}+a_{i+1,j},\ a_{i,j}\}=a_{i,j},
\end{equation}
which is true since $\Delta_{i,n}^{1,j},\ \Delta_{i,n-1}^{2,j},\dots, \Delta_{i,i+1}^{n-i,j}\ge 0$. This completes the inductive step and proves (\ref{1}) and i). \vspace{1mm}

ii) Now conversely, assume that 
$$A=(b_{i,j})_{1\le i,j\le n} \odot p_0\odot  (c_{i,j})_{1\le i,j\le n}\in B_{trop}^{(n)}\odot p_0\odot  B_{trop}^{(n)},$$
\noindent for some $b_{i,j},\ c_{i,j}\in \RR_{\infty}$. We conduct the proof by induction on $n-i$. For the base of induction $i=n$ we note that $b_{n,n}<\infty$, therefore we replace $b_{n,n}$ by $0$ and $c_{1,j},\ 1\le j\le n$ by $c_{1,j}+b_{n,n}$. Thus, $c_{1,j}=a_{n,j}$. As the inductive hypothesis assume that for all $q>i$ it holds
$$\Delta_{q,s}^{k,m}>0,\ q<s\le n,\ 1\le k<m\le n,\ k+s=n+1;$$
$$b_{q,r}=a_{q,n+1-r}-a_{r,n+1-r},\ q\le r\le n;\ c_{n+1-q,j}=a_{q,j},\ n+1-q\le j\le n.$$

At the inductive step we have $a_{i,1}=b_{i,n}+c_{1,1}=b_{i,n}+a_{n,1}$, hence $b_{i,n}=a_{i,1}-a_{n,1}$. After that we have 
$$a_{i,2}=\min\{b_{i,n}+c_{1,2},\ b_{i,n-1}+c_{2,2}\}=\min\{a_{i,1}-a_{n,1}+a_{n,2},\ b_{i,n-1}+a_{n-1,2}\},$$
\noindent provided that $i\le n-2$ (the case $i\ge n-1$ will be considered below as a part of a more general setting). Therefore, it holds $\Delta_{i,n}^{1,2}\ge 0$. If $\Delta_{i,n}^{1,2}= 0$, ii) is justified. Suppose that $\Delta_{i,n}^{1,2}> 0$. Then it holds $b_{i,n-1}=a_{i,2}-a_{n-1,2}$. Continuing in this way by induction on $1\le j\le n-i$, we obtain that 
either

$\bullet$ $\Delta_{i,u}^{n+1-u,r}=0$ for some $1\le r<j,\ 1\le i<u\le n,\ n+1-u<r$ or

$\bullet$ $\Delta_{i,u}^{n+1-u,r}>0$ for all $1\le r<j,\ 1\le i<u\le n,\ n+1-u<r$ and in addition, $b_{i,n+1-r}=a_{i,r}-a_{n+1-r}$.

Then (\ref{3}) implies that 
$$a_{i,j}=\min\{(a_{i,1}-a_{n,1})+a_{n,j},\dots, (a_{i,j-1}-a_{n-j+2,j-1})+a_{n-j+2,j},\ b_{i,n-j+1}+a_{n-j+1,j}\}.$$
Hence it holds $\Delta_{i,n}^{1,j},\ \Delta_{i,n-1}^{2,j},\dots, \Delta_{i,n-j+2}^{j-1,j}\ge 0$. If one of the latter expressions vanishes then ii) is justified. Suppose that $\Delta_{i,n}^{1,j},\ \Delta_{i,n-1}^{2,j},\dots, \Delta_{i,n-j+2}^{j-1,j}>0$. This entails that $b_{i,n-j+1}=a_{i,j}-a_{n-j+1,j}$.

Now consider the case $j=n-i+1$. Then (\ref{3}) implies that
$$a_{i,n-i+1}=\min\{(a_{i,1}-a_{n,1})+a_{n,n-i+1},\dots, (a_{i,n-i}-a_{i+1,n-i})+a_{i+1,j},\ b_{i,i}+c_{n-i+1,n-i+1}\}.$$
\noindent Therefore, it holds $\Delta_{i,n}^{1,n-i+1},\dots, \Delta_{i,i+1}^{n-i,n-i+1}\ge 0$. Again, if one of the latter expressions vanishes then ii) is justified. Suppose that $\Delta_{i,n}^{1,n-i+1},\dots, \Delta_{i,i+1}^{n-i,n-i+1}>0$. Then $a_{i,n-i+1}=b_{i,i}+c_{n-i+1,n-i+1}$. Hence $b_{i,i}<\infty$, so we replace $b_{i,i}$ by $0$ and $c_{n-i+1,j}$ by $c_{n-i+1,j}+b_{i,i},\ n-i+1\le j\le n$. Thus, we obtain that $c_{n-i+1,n-i+1}=a_{i,n-i+1}$.

Finally, consider the case $j\ge n-i+2$. Then (\ref{4}) entails that 
$$a_{i,j}=\min\{(a_{i,1}-a_{i,n})+a_{n,j},\dots, (a_{i,n-i}-a_{i+1,n-i})+a_{i+1,j},\ c_{n-i+1,j}\}.$$
Arguing as above we obtain that either one of $\Delta_{i,n}^{1,j},\dots, \Delta_{i,i+1}^{n-i,j}$ vanishes, or they all are positive, and in this case it holds $c_{n-i+1,j}=a_{i,j}$. This completes the inductive step and the proof of ii).

Observe that to prove  iii) one can follow the proof of ii), taking into account that the case of vanishing one of $\Delta_{i,m}^{n-i+1,j}$ never emerges. $\Box$

\begin{example}
Note that $A\in B_{trop}^{(n)}\odot p_0\odot A\in B_{trop}^{(n)}$ does not necessarily imply (\ref{2}) as demonstrates the following example:
$$A:=\left(\begin{matrix}
0 & -1 & 0 \\ 0 & 0 & 0\\ 
0 & 0 & 0 \end{matrix}\right)
= \left(\begin{matrix} 0 & -1 & 0\\ \infty & 0 & 0\\ \infty & \infty & 0 \end{matrix}\right) \odot p_0 \odot \left(\begin{matrix} 0 & 0 & 0\\ \infty & 0 & 1\\ \infty & \infty & 0\end{matrix}\right)
$$
\noindent in which $\Delta_{1,2}^{2,3}=-1<0$.
\end{example}

\section{Uniqueness of a min-plus Bruhat cell for a tropically non-singular matrix}\label{three}

Recall (see \cite{AGG}) that $n\times n$ matrix $A$ is tropically non-singular if there exists a (unique) permutation $p_A\in S_n$ such that $\sum_{1\le i\le n} a_{i,p_A(i)} < \sum_{1\le i\le n} a_{i,p(i)}$ for any permutation $p\neq p_A$.

\begin{proposition}
For a tropically non-singular matrix $A$ if $A\in B_{trop}^{(n)}\odot p\odot  B_{trop}^{(n)},\ p\in S_n$ then $p=p_A$. 
\end{proposition}

{\bf Proof}. We use the following statement: if $A=U\odot V$ for some $n\times n$ matrices $U,\ V$ then $U,\ V$ are also tropically non-singular and $p_A=p_U\circ p_V$ (see Proposition 3.4 \cite{M}, implicitly one can also find it in \cite{B}). Therefore $p_A=p$. $\Box$ \vspace{2mm}

For a classical (not necessarily non-singular) $n\times n$ matrix $D$ there exists a (unique) permutation $p_D$ such that $D\in B^{(n)}p_DB^{(n)}$, and if $D\in B^{(n)}pB^{(n)}$ then $p_D \preceq p$  \cite{G81}, \cite{G82}. Is there a tropical analogue of this assertion?

\end{document}